	\magnification=\magstep1
\def\blackbox{\hbox{\vrule width6pt height7pt depth1pt}}
\def\qed{\hfill\blackbox}
\def\precsim{\lower.5ex\hbox{$ 
	\ \mathop{\buildrel \prec\over{\scriptstyle\sim}}\nolimits\ $}}
\def\rf#1{{\rm [#1]}}
\def\IN{\mathop{{\rm I}\kern-.2em{\rm N}}\nolimits}
\def\proof{\noindent {\it Proof.}\enspace}
\def \wtil{\widetilde}
\topinsert\vskip.5in\endinsert
\centerline{\bf On quotients of Banach spaces having}
\smallskip
\centerline{\bf shrinking unconditional bases}
\bigskip
\centerline{by E. Odell\footnote*
{Research partially supported by NSF Grant DMS-8903197.}}

\vskip.5in
{\narrower\smallskip
\centerline{\bf Abstract}
\medskip

It is proved that if a Banach space $Y$ is a quotient of a Banach space having
a shrinking unconditional basis, then every normalized weakly null sequence in 
$Y$ has an unconditional subsequence.  The proof yields the corollary that
every quotient of Schreier's space is $c_o$-saturated.
\smallskip}
\bigskip

\noindent{\bf \S 0. Introduction}.

We shall say that a Banach space $Y$ has {\it property (WU)\/} if every 
normalized weakly null sequence in $Y$ has an unconditional subsequence.
The well known example of Maurey and Rosenthal \rf{MR} shows that not every 
Banach space has property (WU) (see also \rf{O}).  W.B.~Johnson \rf{J}
proved that if $Y$ is a quotient of a Banach space $X$ having a shrinking
unconditional f.d.d.\ and the quotient map does not fix a copy of $c_0$,
then $Y$ has (WU).  Our main result extends this (and solves Problem IV.1
of \rf{J}).

\proclaim Theorem A.  Let $X$ be a Banach space having a shrinking
unconditional finite dimensional decomposition.  Then every quotient of $X$
has property (WU).
\par

Of course such an $X$ will itself have property (WU).  
Furthermore, if $(E_n)$ is an
unconditional f.d.d.\ (finite dimensional decomposition) for $X$, then $(E_n)$
is shrinking if and only if $X$ does not contain $\ell_1$.

The proof of Theorem A yields 

\proclaim Theorem B.  Let $Y$ be a Banach space which is a quotient of $S$, 
the Schreier space.  Then $Y$ is $c_o$-saturated.

$Y$ is said to be {\it $c_o$-saturated\/} if every infinite dimensional 
subspace of $Y$ contains an isomorph of $c_0$.

Our notation is standard as may be found in the books of Lindenstrauss and
Tzafriri \rf{LT 1,2}.  The proof of Theorem A is given in \S1 and the proof 
of Theorem B appears in \S2.  \S3 contains some open problems.  We thank 
H.~Rosenthal and T.~Schlumprecht for useful conversations regarding the 
results contained herein.
\bigskip

\noindent{\bf \S 1. The proof of Theorem A}.

Let $T$ be a bounded linear operator from $X$ onto $Y$ where $X$ has a
shrinking unconditional f.d.d., $(\buildrel \approx \over E_i)$.  By 
renorming if necessary we may suppose that 
$(\buildrel \approx \over E_i)$ is 1-unconditional.
$Y^*$ is separable and so by a theorem of Zippin \rf{Z} we may assume
that $Y$ is a subspace of a Banach space $Z$ possessing a bimonotone
shrinking basis, $(z_i)$.  Fix $C >0$ such that
$$T(C B_a X)\supseteq B_a Y \equiv \{y\in Y: \|y\|\le 1\}\ .$$

Recall that $(\wtil E_i)$ is a {\it blocking\/} of $(\buildrel \approx \over
E_i)$ if there exist integers $0 = q_0 <q_1 <q_2 <\cdots$ such that
$\wtil E_i = [\buildrel \approx \over E_j]^{q_i}_{j=q_{i-1}+1}$ for all $i$
(where $[\cdots]$ denotes the closed linear span).
Similarly, $\wtil F_i = [z_j]^{q_i}_{j=q_{i-1}+1}$ defines a blocking of
$(z_i)$.

Fix a sequence $\varepsilon_{-1} > \varepsilon_0 > \varepsilon_1 >
\varepsilon_2 > \cdots$ converging to 0 which satisfies
$$\sum_{i=-1}^\infty \varepsilon_i < 1/4 \ \hbox{ and }\  \sum^\infty_{i=p}
(4i+2) \varepsilon_i < \varepsilon_{p-1}\ \hbox{ for }\  p \ge 0\ .\leqno(1.1)$$

Then choose $\tilde \varepsilon_0 > \tilde \varepsilon_1 > \cdots$ converging
to 0 which satisfies
$$4p\tilde\varepsilon_p< \varepsilon_{p+2}\ \hbox{ for }\ p\ge1\ \hbox{ and }
\ \sum^\infty_{j=p+1}\tilde\varepsilon_j <\tilde \varepsilon_p\ \hbox{ for }
\ p \ge 0\ .\leqno(1.2)$$

Our first step is the blocking technique of Johnson and Zippin.

\proclaim Lemma 1.1 (\rf{JZ 1,2}).  There exist blockings $(\wtil E_i)$
and $(\wtil F_i)$ of $(\buildrel \approx \over E_i)$ and $(z_i)$,
respectively, such that if $(\wtil Q_i)$ is the sequence of finite rank
projections on $Z$ associated with $(\wtil F_i)$ then
$$\qquad \hbox{For all } i\in \IN \hbox{ and } x\in \wtil E_i \hbox{ with }
\Vert x\Vert \le C, \hbox{ we have } \Vert \wtil Q_j Tx\Vert < \tilde 
\varepsilon_{\max(i,j)} \hbox{ if } j\ne i, i-1\ .\leqno(1.3)$$

Roughly, this says that $T \wtil E_i$ is essentially contained in
$\wtil F_{i-1} +\wtil F_i$ (where $\wtil F_0 = \{0\}$).  Let $(y''_i)$
be a normalized weakly null sequence in $Y$.  Choose a subsequence $(y''_i)$
of $(y_i)$ and a blocking $(F_i)$ of $(\wtil F_i)$, given by
$F_i = [\wtil F_j]^{q_i}_{j=q_{i-1}+1}$, such that if $Q_i = \sum^{q_i}
_{j=q_{i-1}+1} \wtil Q_j$ is the sequence of finite rank projections on
$Z$ associated with $(F_i)$, then
$$\Vert Q_j y'_i\Vert < \tilde \varepsilon_{\max(i,j)} \hbox{ if } i\ne j\ .
\leqno(1.4)$$
Roughly, $y'_i$ is essentially in $F_i$.  Furthermore we may assume that
$$\Vert \sum a_i y'_i\Vert = 1 \hbox{ implies } \max \vert a_i\vert \le 2\ .
\leqno(1.5)$$

Let $(E_i)$ be the blocking of $(\wtil E_i)$ given by the {\it same}
sequence $(q_i)$ which defined $(F_i)$, $E_i = [\wtil E_j]^{q_i}_{j=q_{i-1}+1}
$.

We begin with a sequence of elementary technical yet necessary lemmas.

For $I\subseteq \IN$ we define $Q_I = \sum_{j\in I} Q_j$ and set
$Q_\phi =0$.

\proclaim Lemma 1.2.  Let $0 \le n < m$ be integers and let $y =\sum
_{i\notin (n,m)} a_i y'_i$ with $\Vert y\Vert =1$.  Then for $j\in
(n,m)$, $\Vert Q_j y\Vert < \varepsilon_j$ and $\Vert Q_{(n,m)} y\Vert 
< \varepsilon_n$.
\par

\proof  Let $n <j < m$.  Then by (1.5), (1.4), (1.2) and (1.3), 
$$\eqalign{\Vert Q_j y
\Vert &\le 2\Bigl(\sum_{i\le n} \Vert Q_j y'_i\Vert +\sum_{i\ge m} \Vert 
Q_j y'_i\Vert\Bigr)\cr
&< 2(n \tilde \varepsilon_j +\tilde \varepsilon_{m-1})\cr
&\le (2j +2)\tilde \varepsilon_j \le 4j \tilde \varepsilon_j < \varepsilon_j
\cr}\ .
$$
Thus $\Vert Q_{(n,m)} y\Vert < \sum_{j\in(n,m)}\varepsilon_j <
\varepsilon_n$ by (1.1).~\qed

\proclaim Lemma 1.3.  Let $0 = p_0 < r_0 =1 < p_1 < r_1 < p_2 < r_2 < \cdots$
be integers and let $y = \sum^\infty_{i=1} a_i y'_{p_i}$ 
with $\Vert y\Vert =1$.  Then
for $i\in \IN$,
$$\Vert Q_{[r_{i-1},r_i)} y - a_i y'_{p_i}\Vert < \varepsilon_{p_{i-1}-1}\ .$$
\par

\proof  
$$\eqalignno{
&\Vert Q_{[r_{i-1},r_i)} y - a_i y'_{p_i}\Vert\cr
&\le \Vert Q_{[r_{i-1},r_i)}\sum_{j\ne i} a_j y'_{p_j}\Vert + \Vert Q_{[r_{i-1},
r_i)} a_i y'_{p_i} - a_i y'_{p_i}\Vert\cr
&\hbox{ which by lemma 1.2 is }\cr
&< \varepsilon_{r_{i-1}-1} + \Vert Q_{[1,r_{i-1})} a_i y'_{p_i}\Vert +\Vert
Q_{[r_i,\infty)} a_i y'_{p_i}\Vert\cr
&< \varepsilon_{r_{i-1}-1} + 2\sum_{k <r_{i-1}} \Vert Q_k y'_{p_i}\Vert +
2\varepsilon_{r_i-1} \hbox{ (by (1.5) and lemma 1.2) }\cr
&< \varepsilon_{r_{i-1}-1} +2(r_{i-1}-1) \tilde \varepsilon_{p_i} +
2\varepsilon_{r_i-1} \hbox{ (by (1.4))}\cr
&\le \varepsilon_{p_{i-1}} +2p_i \tilde \varepsilon_{p_i} 
+2\varepsilon_{p_i} < \varepsilon_{p_
{i-1}} +4 \varepsilon_{p_i} \hbox{ (by (1.2))}\cr
&<\varepsilon_{p_{i-1} -1} \hbox{ (by 1.1)}\ .&\blackbox\cr}$$

\proclaim Lemma 1.4.  Let $i\in \IN$, $x\in E_i$ and $\Vert x\Vert \le C$.  Then
$$\eqalign{
\Vert Q_j Tx\Vert &< \varepsilon_{\max(i,j)} \hbox{ if } j \ne i,i-1\ ,\cr
\Vert Q_{[1,i-2]} Tx\Vert &< \varepsilon_{i-1} \hbox{ and } \Vert Q_{[i,\infty)}
Tx\Vert <\varepsilon_{i-1}\ .\cr}$$
\par

\proof Let $x = \sum_{\ell \in (q_{i-1},q_i]} \omega_{\ell}$ with
$\omega_{\ell} \in \wtil E_{\ell}$.  
$$\eqalign{
\Vert Q_j Tx\Vert &\le \sum_{k\in(q_{j-1},q_j]}\sum_{\ell \in (q_{i-1}, q_i]}
\Vert \wtil Q_k T\omega_\ell \Vert\cr
\hbox{(if $j <i-1$)}: &< \sum_{k\in (q_{j-1}, q_j]} \sum_{\ell \in (q_{i-1},
q_i]}\tilde \varepsilon_\ell \hbox{ (by (1.4))}\cr
&< q_j \tilde \varepsilon_{q_{i-1}} < q_{i-1} \tilde \varepsilon_{q_{i-1}}
< \varepsilon_{q_{i-1} +2} < \varepsilon_i
\hbox{ using\ (1.2) }\cr
\hbox{ (if $j >i$)}:  &< \sum_{k\in(q_{j-1}, q_j]}\sum_{\ell \in (q_{i-1}, q_i]}
\tilde \varepsilon_k\cr
&< \sum_{k\in(q_{j-1},q_j]}q_i \tilde \varepsilon_k \le q_i \tilde \varepsilon
_{q_{j-1}}\cr
&\le q_{j-1} \tilde \varepsilon_{q_{j-1}} < \varepsilon_{q_{j-1} +2} \le
\varepsilon_{j+1} < \varepsilon_j\ .\cr}$$

Finally, 
$$\Vert Q_{[1,i-2]} Tx\Vert \le \sum^{i-2}_{k=1} \Vert Q_k Tx\Vert < 
\sum^{i-2}_{k=1} \varepsilon_i = (i-2)\varepsilon_i < \varepsilon_{i-1}$$
and
$$\Vert Q_{[i,\infty)} Tx\Vert \le \sum^\infty_{k=i} \Vert Q_k Tx\Vert
<\sum^\infty_{k=i} \varepsilon_k < \varepsilon_{i-1}\ .\eqno\blackbox$$
\medskip

\proclaim Lemma 1.5.  Let $\Vert x\Vert \le C$, $x = \sum_{k\ne j,j+1}
\omega_k$ where $\omega_k \in E_k$ for all $k$.  Then
$$\Vert Q_j Tx\Vert < \varepsilon_{j-1}\ .$$
\par

\proof By lemma 1.4, 
$$\eqalignno{
\Vert Q_j Tx\Vert &\le \sum_{k\ne j, j+1}\Vert Q_j T\omega_k\Vert < \sum_{k<j}
\varepsilon_j +\sum_{k>j+1} \varepsilon_k\cr
&< (j-1) \varepsilon_j +\varepsilon_j = j\varepsilon_j < 
\varepsilon_{j-1}\ .&\blackbox\cr}$$

\proclaim Lemma 1.6.  Let $1\le n < m$ and $x = \sum \omega_j$, $\Vert x\Vert
\le C$, with $\omega_j \in E_j$ for all $j$.  Suppose that $\Vert Q_j Tx\Vert
< 2 \varepsilon_{j-1}$ for $n <j <m$.  Let $a_{j-1} = Q_{j-1} T\omega_j$ and
$b_j = Q_j T\omega_j$.  Then
\itemitem{a)}  $\Vert a_j +b_j\Vert < 3 \varepsilon_{j-1}$ for $n < j<m$ and
\medskip
\itemitem{b)}  $\Vert \sum_{j\in(r,s]} T\omega_j -(a_r +b_s)\Vert
< 5 \varepsilon_{r-1}$ if $n < r < s < m$.
\par

\proof  a) Let $n <j <m$.  By lemma 1.5, $\Vert Q_j Tx -(a_j +b_j)\Vert =
\Vert Q_j(\sum_{i\ne j,j+1} T\omega_i)\Vert < \varepsilon_{j-1}$.
\itemitem{}  Since $\Vert Q_j Tx\Vert < 2\varepsilon_{j-1}$, a) follows.
\smallskip
\itemitem{b)}  Let $n <r <s <m$ and let $j\in (r,s]$.  Then $T\omega_j =
a_{j-1} +b_j +\gamma_j$ where $\Vert \gamma_j\Vert < 2\varepsilon_{j-1}$
by lemma 1.4.
Thus 
$$\eqalignno{
&\Vert \sum^s_{r+1} T\omega_j -(a_r +b_s)\Vert\cr
&\le \Vert a_r +b_{r+1} + a_{r+1} + b_{r+2} +\cdots + a_{s-1}+ b_s -
(a_r +b_s)\Vert +\sum^s_{j=r+1} 2\varepsilon_{j-1}\cr
&< \sum^{s-1}_{r+1} \Vert a_j +b_j\Vert +2 \varepsilon_{r-1}\cr
&< \sum^{s-1}_{r+1} 3\varepsilon_{j-1} +2\varepsilon_{r-1} \hbox{ (by a))}\cr
&< 5 \varepsilon_{r-1}\ .&\blackbox\cr}$$

We next come to the key lemma.  Let $(P_j)$ be the sequence of finite rank
projections on $X$ associated with $(E_j)$.  For $I\subseteq \IN$, we let
$P_I = \sum_{i\in I} P_i$.  
\medskip
\item{}  {\it notation:}  If $x =\sum x_j\in X$ with $x_j\in E_j$ for all $j$
and $\overline x \in X$, we define
$$\overline x \precsim x \hbox{ if } \overline x = \sum a_j x_j \hbox{ with }
0 \le a_j \le 1 \hbox{ for all } j\ .$$

\proclaim Lemma 1.7.  Let $n\in \IN$ and let $\varepsilon > 0$.  There exists
$m\in \IN$, $m > n+1$, such that whenever $x\in C Ba X$ with $\Vert Q_j
Tx\Vert < 2\varepsilon_{j-1}$ for all $j\in (n,m)$ then:  there exists
$\overline x\precsim x$ with
\itemitem{1)}  $\Vert Tx -T\overline x\Vert < \varepsilon$ and
\smallskip
\itemitem{2)}  $P_r\overline x =0$ for some $r\in (n,m)$.
\par

\noindent{\it Remark}.  Lemma 1.7 is the main difference between our result
and Johnson's earlier special case \rf{J}.  In the case where $T$ does not fix 
a copy of $c_0$, Johnson showed that one could take $\overline x =
x -P_r(x)$ for some $r\in (n,m)$.

The proof of lemma 1.7 requires the following key

\proclaim Sublemma 1.8.  Let $n\in \IN$ and $\varepsilon >0$.  There exists
an integer $m = m(n,\varepsilon) > n+1$ satisfying the following.  Let
$x\in C Ba X$, $x = \sum \omega_j$ with $\omega_j\in E_j$ for all $j$.
Assume in addition that $\Vert Q_j Tx\Vert < 2 \varepsilon_{j-1}$ for 
$j\in (n,m)$ and set $a_{j-1} = Q_{j-1} T\omega_j$ and $b_j = Q_j T\omega_j$.
Then there exist $k\in \IN$ and integers $n < i_1 < \cdots < i_k < m$ such
that
$$k^{-1}\Vert a_{i_1} +a_{i_2} +\cdots +a_{i_k}\Vert <\varepsilon\ .\leqno
(1.6)$$
\par

\noindent{\it Proof of Lemma 1.7}.  Let $n\in \IN$ and $\varepsilon >0$.
Choose $n_0 \ge n$ such that
$$\varepsilon_{n_0} < \varepsilon/12\ .\leqno(1.7)$$ 
Let $m_1 = m(n_0+1, \varepsilon/3)$ be given by the sublemma and let
$m = m(m_1, \varepsilon/3)$.

Let $x = \sum \omega_j\in C Ba X$ with $\omega_j\in E_j$ for all $j$ and
suppose that $\Vert Q_j Tx\Vert < 2\varepsilon_{j-1}$, $a_{j-1} = Q_{j-1}
T\omega_j$ and $b_j = Q_j T\omega_j$ for $j\in (n,m)$.  By our choice of
$m$ there exist integers $k$ and $K$ and integers $n\le n_0 < n_0 +1 < i_1
< i_2 <\cdots < i_k <m_1 <j_1 <\cdots < j_K <m$ such that
$$\leqalignno{
&k^{-1}\Vert a_{i_1} + \cdots + a_{i_k} \Vert < \varepsilon/3 \hbox{ and}
&(1.8)\cr
&K^{-1}\Vert a_{j_1} +\cdots + a_{j_K}\Vert < \varepsilon/3\ .&(1.9)\cr}$$
Define 
$$\eqalign{
\overline x &= \sum^{i_1}_1 \omega_j +{k-1\over k} \sum^{i_2}_{i_1+1}
\omega_j + \cdots +{1\over k}\sum^{i_k}_{i_{k-1}+1}
+{0\over k}
\sum^{j_1}_{i_k+1} \omega_j\cr 
&\qquad +{1\over K} \sum^{j_2}_{j_1+1} \omega_j +\cdots +
{K\over K} \sum^\infty_{j_k+1}\omega_j\ .\cr}$$
Clearly (2) holds and we are left to check (1).
$$\eqalign{
&\Vert Tx -T\overline x\Vert = \Vert {1\over k} \sum^{i_2}_{i_1+1}T\omega_j
+{2\over k} \sum^{i_3}_{i_2+1} T\omega_j\cr
&\qquad +\cdots +{k\over k}\sum^{j_1}_{i_k+1} T\omega_j +{K-1\over K} \sum^{j_2}
_{j_1+1}T\omega_j
+ \cdots + {1\over K} \sum^{j_K}_{j_{K-1}+1}
T\omega_j \Vert\ .\cr}$$

Thus by lemma 1.6,
$$\eqalign{
&\Vert Tx -T\overline x\Vert \le \Vert {1\over k} a_{i_1} +{1\over k} b_{i_2} 
+{2\over k} a_{i_2} +{2\over k} b_{i_2}\cr
&\qquad + \cdots + {k\over k} a_{i_k} +{K\over K} b_{j_1} +{K-1\over K} a_{j_1}
+{K-1\over K} b_{j_2}
+ \cdots + {1\over K} a_{j_{K-1}} +{1\over K} b_{j_K}\Vert\cr
&\qquad + k^{-1} \sum_{j=1}^k 5j \varepsilon_{i_{j-1}} +K^{-1} \sum^K_{\ell =1}
5\ell \varepsilon_{j_\ell -1}\ .\cr}$$

Now ${k^{-1} \sum^k_{j=1} 5j \varepsilon_{i_j-1} \le 5 \sum^k_{j=1} 
\varepsilon_{i_j-1} < \varepsilon_{i_1-2} \le \varepsilon_{n_0}}$ and
$K^{-1} \sum^K_{\ell =1} 5\ell \varepsilon_{j_\ell-1} < \varepsilon_{n_0}$
as well.

Thus
$$\eqalign{
\Vert Tx -T\overline x\Vert &< k^{-1} \Vert a_{i_1} +\cdots + a_{i_k} \Vert +
K^{-1}\Vert b_{j_1} +\cdots + b_{j_K}\Vert\cr
&\qquad +k^{-1} \sum^k_{j=2} \Vert b_{i_j} +a_{i_j}\Vert + K^{-1}
\sum^{K-1}_{\ell =1} \Vert b_{j_\ell} +a_{j_\ell}\Vert 
+2\varepsilon_{n_0}\ .\cr}$$

From (1.8), (1.9) and lemma 1.6 we obtain
$$\eqalign{
\Vert Tx -T\overline x\Vert &< {\varepsilon\over 3} +{\varepsilon\over 3} +
\sum^k_{j=2} 3\varepsilon_{i_j-1} 
+ \sum^{K-1}_{\ell =1} 3\varepsilon_{j_\ell-1} +2\varepsilon_{n_0}\cr
&< {2\varepsilon\over 3} + \varepsilon_{n_0} + \varepsilon_{n_0}
+ 2\varepsilon_{n_0} <\varepsilon\cr}$$
(by (1.7)).\qed
\medskip
\noindent{\it Proof of Sublemma 1.8}.  If the sublemma fails then by a
standard compactness argument we obtain $\omega_j\in E_j$ for $j\in \IN$ such
that for all $m$,
$$\Vert \sum^m_{j=1} \omega_j\Vert \le C \hbox{ and } \Vert Q_j T(\sum^m_{i=1}
\omega_i)\Vert \le 3\varepsilon_{j-1}$$
if $n <j <m$.  The extra $\varepsilon_{j-1}$ comes from an application of
lemma 1.5.  Furthermore setting $Q_{j-1} T\omega_j = a_j$ and $Q_j T\omega_j
=b_j$ for $j\in \IN$, then for all $k$ and all $n <i_1 <\cdots <i_k$ we have
$$k^{-1} \Vert a_{i_1} +\cdots + a_{i_k}\Vert \ge \varepsilon\ .\leqno(1.10)$$

Now $a_j \in F_j$ and $(F_j)$ is a shrinking f.d.d.  Thus $(a_j)_{j>n}$ is
a seminormalized weakly null sequence.  By (1.10) any spreading model of a
subsequence of $(a_j)$ must be equivalent to the unit vector basis of
$\ell_1$ (see \rf{BL} for basic information on spreading models).  In
particular we can choose an even integer $k$ and integers
$n <i_1 <\cdots <i_k$ such that
$$\Vert a_{i_1} - a_{i_2} +\cdots + a_{i_{k-1}} -a_{i_k}\Vert > C\Vert T\Vert 
+1\ .\leqno(1.11)$$ 
However,
$$\eqalign{
C\Vert T\Vert &\ge \Vert T(\sum^{i_2}_{i_1+1} \omega_j +\sum^{i_4}_{i_3+1}
+\cdots +\sum^{i_k}_{i_{k-1}+1}\omega_j)\Vert\cr
&\ge \Vert a_{i_1} +b_{i_2}+a_{i_3} +b_{i_4} 
+\cdots + a_{i_{k-1}} +b_{i_k}\Vert\cr
&\qquad -5 \sum^k_{j=1} \varepsilon_{i_j-1} \hbox{ (by lemma 1.6)}\ .\cr}$$
Now $5 \sum^k_{j=1} \varepsilon_{i_j-1} < \varepsilon_{i_1-2}$
and by lemma 1.6 and (1.11)
$$\eqalign{
\Vert a_{i_1} &+b_{i_2} +\cdots + a_{i_{k-1}} +b_{i_k}\Vert\cr
&\ge \Vert a_{i_1} -a_{i_2} +a_{i_3} -a_{i_4} +\cdots + a_{i_{k-1}}-a_{i_k}
\Vert\cr
&-\sum^{k/2}_{j=1} \Vert a_{i_{2j}}+b_{i_{2j}}\Vert > C\Vert T\Vert +1 -
\sum^{k/2}_{j=1} 3\varepsilon_{i_{2j}-1}\ .\cr}$$
Thus
$$\eqalign{
C\Vert T\Vert &> C\Vert T\Vert +1 -\varepsilon_{i_1-2}-\varepsilon_{i_2-2}\cr
&\ge C\Vert T\Vert +1 -2\varepsilon_{i_1-2} > C\Vert T\Vert\ ,\cr}$$
which is impossible.\qed
\bigskip

\noindent{\it Completion of the proof of Theorem A}.  

Let the integer $m$ given by lemma 1.7 be denoted by $m = m(n;\varepsilon)$.
Choose $1 <p_1 <p_2 <\cdots$ such that for all $i$, $p_{i+1} -1 \ge m
(p_i; \varepsilon_{p_i})$.  Let $(y_i) = (y'_{p_i})$.  We shall prove that
$(y_i)$ is unconditional.

Let $y = \sum a_i y_i, \Vert y\Vert =1$, $x\in C Ba X$, $Tx =y$ and let
$x = \sum^\infty_{i=0} g_i$ where $g_0 = P_{[1,p_1)}x$ and $g_i =
P_{[p_i, p_{i+1})}x$ for $i \ge 1$.  We shall apply lemma 1.7 to each
$g_i$ for $i\ge 1$.  Fix $i \ge 1$ and let $(n,m) = (p_i, p_{i+1} -1)$.
Let $j\in (n,m)$.  Then $\Vert Q_jy\Vert <\varepsilon_j$ by lemma 1.2.
Thus $\Vert Q_j Tx\Vert = \Vert Q_j Tg_i +Q_jT \sum_{k\ne i} g_k\Vert
<\varepsilon_j$.
However $\Vert Q_j T \sum_{k\ne i} g_k\Vert < \varepsilon_{j-1}$ by
lemma 1.5 so $\Vert Q_j T g_i\Vert < \varepsilon_{j-1} +\varepsilon_j
< 2\varepsilon_{j-1}$.  Thus by lemma 1.7 there exist $\overline g_i 
\precsim g_i$
and $r_i\in (p_i, p_{i+1}-1)$ such that $P_{r_i} \overline g_i =0$ and
$\Vert T g_i - T\overline g_i\Vert < \varepsilon_{p_i}$ for all $i\in \IN$.

Let $\overline x = \sum^\infty_{i=0} \overline g_i = \sum^\infty_{i=1}
\overline x_i$ where $\overline g_0 = g_0$ and $\overline x_i = p_{[r_{i-1},
r_i)}\overline x$ for $i\in \IN\ (r_0 =1)$.  Of course, $\overline x_i = 
P_{(r_{i-1},r_i)}\overline x$ if $i >1$.
\medskip

\noindent{\it Claim\/}:  $\Vert T\overline x_i - a_i y_i\Vert < 4 \varepsilon_
{p_{i-1}-1}$ for $i\in \IN$.

Indeed $\Vert Q_{[r_{i-1},r_i)} y-a_i y_i\Vert < \varepsilon_{p_{i-1}-1}$
by lemma 1.3.  Thus the claim follows from the 
\medskip

\noindent{\it Subclaim\/}:  $\Vert Q_{[r_{i-1},r_i)} Tx -T\overline x_i\Vert <
3\varepsilon_{p_{i-1}-1}$.

To see this we first note that 
$$\eqalign{
\Vert Q_{[r_{i-1},r_i)} Tx &-Q_{[r_{i-1},r_i)} T(g_{i-1} +g_i +g_{i+1})\Vert\cr
&\le \sum_{k\in [r_{i-1},r_i)} \Vert Q_k \sum_{j\ne i-1,i,i+1} T g_j\Vert\cr
&< \sum_{k\in [r_{i-1},r_i)} \varepsilon_{k-1} \qquad \hbox{(by lemma 1.5)}\cr
&< \varepsilon_{r_{i-1}-1}\ .\cr}$$
Also
$$\eqalign{
\Vert Q_{[r_{i-1},r_i)} &T(g_{i-1} +g_i +g_{i+1}) -Q_{[r_{i-1},r_i)} 
T(\overline g_{i-1} +
\overline g_i +\overline g_{i+1})\Vert\cr
&\le \Vert T(g_{i-1} +g_i +g_{i+1} -\overline g_{i-1} - \overline g_i -
\overline g_{i+1})\Vert\cr
&< \varepsilon_{p_{i-1}} + \varepsilon_{p_i} + \varepsilon_{p_{i+1}} < 
\varepsilon_{p_{i-1}-1}\ .\cr}$$
Finally, applying lemma 1.5 again we have
$$\eqalign{
\Vert Q_{[r_{i-1},r_i)} &[T(\overline g_{i-1} +\overline g_i +\overline g_{i+1})
-T(\overline x_i)]\Vert\cr
&< \varepsilon_{r_{i-1}-1}\ , \hbox{ and the subclaim follows }\ .\cr}$$

Let $\delta_i = \pm 1$.  Then
$$\eqalignno{
\Vert \sum \delta_i a_i y_i\Vert &\le \Vert \sum \delta_i(a_i y_i -T\overline
x_i)\Vert\cr
&+ \Vert \sum \delta_i T\overline x_i\Vert\cr
&< \sum 4\varepsilon_{p_{i-1}-1} + \Vert T\Vert\ \Vert \sum \delta_i \overline
x_i\Vert\cr
&\qquad \hbox{(by the claim)}\cr
&\le 1 + C\Vert T\Vert\ .&\blackbox\cr}$$

The proof of Theorem A yields the following

\proclaim Proposition 1.9.  Let $X$ have a shrinking $K$-unconditional f.d.d.
$(E_i)$ and let $T$ be a bounded linear operator from $X$ onto $Y$.  Let $T(C BaX)
\supseteq Ba Y$.  Then if $\varepsilon_i \downarrow 0$ and if $(y'_i)$ is a
normalized weakly null basic sequence in $Y$ there exists a subsequence 
$(y_i)$ of $(y'_i)$ and integers $p_1 <p_2 <\cdots$ with the following
property.  Let $\Vert \sum a_i y_i\Vert \le 2$.  Then there exists 
$x = \sum x_i \in 2C KBa X$, $(x_i)$ a block basis of $(E_i)$, such that
$$\Vert Tx_i -a_i y_i\Vert < \varepsilon_i \quad \hbox{ for all } i\ .$$
Moreover there exist $(r_i)$ with $0 = r_0 < p_1 < r_1 < p_2 < r_2 < \cdots$
such that $x_i \in [E_j]_{j\in(r_{i-1},r_i)}$ for all $i$.
\par

\proclaim Corollary 1.10.  Let $X$ have a shrinking $K$-unconditional f.d.d.
and let $T$ be a bounded linear operator from $X$ onto the Banach space $Y$.
Then $Y$ contains $c_0$ if and only if $T$ fixes a copy of $c_0$.
\par

\proof  If $Y$ contains $c_o$ then there exists (see \rf{Ja}) 
$(y_i)$, a normalized sequence in $Y$, with $2^{-1} \le \Vert
\sum a_i y_i\Vert \le 2$ if $(a_i)\in S_{c_0}$, the unit sphere of $c_o$. 
Let $\varepsilon_i \downarrow
0$ with $\sum \varepsilon_i < 1$.  We may assume that $(y_i)$ satisfies the
conclusion of proposition 1.9.  Thus for all $n\in \IN$ there exist
$$0 = r^n_0 < p_1 < r^n_1 < p_2 < r^n_2 < \cdots$$ 
and $x^n_i\in [E_j]_{j\in(r^n_{i-1},r^n_i)}$ such that if $x^n = \sum_
{i\le n} x^n_i$, then $\Vert x^n\Vert \le 2CK$ and $\Vert Tx^n_i - y_i\Vert
< \varepsilon_i$ for $i\le n$.

By passing to a subsequence $(x^{n_k})$ we may assume $\lim_{k\to \infty}
r^{n_k}_i = r_i$ and $\lim_{k\to \infty} x^{n_k}_i = x_i$ exist for
all $i\in \IN$.  Thus $x_i\in [F_j]_{j\in (r_{i-1},r_i)}$ 
with $r_0 = 0 <r_1 <r_2 <\cdots$,
$\Vert Tx_i - y_i\Vert <\varepsilon_i$ for all $i$ and ${\sup_n}\Vert
\sum^n_1 x_i\Vert <\infty$.  It follows that $(x_i)$ is equivalent to the
unit vector basis of $c_0$.  Moreover if we choose $\omega_i \in \varepsilon_
i C Ba X$ with $T\omega_i = y_i -Tx_i$ then $T(x_i +\omega_i) = y_i$ and
some subsequence of $(x_i +\omega_i)$ is also a $c_0$ basis.  Hence $T$
fixes $c_0$.\qed
\bigskip

\noindent{\bf \S 2.  The proof of Theorem B}.
\smallskip

We begin by recalling the definition of the Schreier space $S$ \rf{S}.  Let
$c_{00}$ be the linear space of all finitely supported real valued sequences.
For $x =(x_i)\in c_{00}$ set
$$\Vert x\Vert = \max_i \{\sum^p_{i=1} 
\vert x_{k_i}\vert: p\in \IN \hbox{ and }
p \le k_1 <\cdots < k_p\}\ .$$
$S$ is the completion of $(c_{00}, \Vert \cdot \Vert)$.  We let $\Vert x
\Vert_0$ denote the $c_0$-norm of $x$.  The unit vector basis $(e_n)$ is
a shrinking 1-unconditional basis for $S$.  $S$ can be embedded into
$C(\omega^\omega)$ and thus $S$ is $c_0$-saturated.

Theorem B will follow from a quantitative version, Theorem B$'$ (below).
Given a sequence $(x_n)$, $\lambda > 0$ and $F$ a finite nonempty subset
of $\IN$, $y = \lambda \sum_{n\in F} x_n$ is said to be a {\it 1-average}
of $(x_n)$.  We say that a Banach space $X$ has property-$S(1)$ if every
normalized weakly null sequence in $X$ admits a block basis of 1-averages 
which is equivalent to the unit vector basis of $c_0$.  $S$ has property-$S(1)$.

\proclaim Theorem B$'$.  Let $Y$ be a quotient of $S$.  Then $Y$ has
property-$S(1)$.
\par

We shall use the following simple 
\proclaim Lemma 2.1.  Let $(x_n)$ be a normalized weakly null sequence in
$S$ with $\lim_n\Vert x_n\Vert_0 =0$.  Then some subsequence of $(x_n)$
is equivalent to the unit vector basis of $c_0$.
\par

Let $T$ be a bounded linear operator from $S$ onto a Banach space $Y$ and
let $(y'_i)$ be a normalized weakly null basic sequence in $Y$.  
Let $T(C Ba S) \supseteq Ba Y$.

\proclaim Lemma 2.2.  If no block basis of 1-averages of $(y'_i)$ is 
equivalent to the unit vector basis of $c_0$, 
then there exists
$\delta > 0$ such that if $x\in 3C Ba S$, $Tx$ is a 1-average of 
$(y'_i)$ and $\Vert Tx
\Vert > 1/3$ then $\Vert x\Vert_0 > \delta$.
\par

\proof If no such $\delta$ exists then there exists $(x_i)\subseteq 3C Ba S$
with $\lim_i \Vert x_i\Vert_0 =0$, $\Vert Tx_i\Vert > {1\over 3}$
and $Tx_i$ a 1-average of $(y'_i)$
for all $i$.  By lemma 2.1 there exists a subsequence $(x'_i)$ of $(x_i)$
which is equivalent to the unit vector basis of $c_0$.  By passing to a
further subsequence we may assume that $(Tx'_i)$ is a seminormalized weakly
null basic sequence in $[(y'_i)]$.  Thus $(Tx'_i)$ is also equivalent to
the unit vector basis of $c_0$.\qed
\medskip

\noindent{\it Proof of Theorem $B'$}.  Let $(y'_i)$ be a normalized weakly null
sequence in $Y$.  If $(y'_i)$ fails the $S(1)$ property, choose $\delta
>0$ by lemma 2.1.  Let $(\varepsilon_i)^\infty_{i=1}$ be a sequence of
positive numbers satisfying
$$\sum^\infty_{i=1} \varepsilon_i < \min(\delta/(2C), 1)\ .\leqno(2.1)$$

Let $(y_i)$ be the subsequence of $(y'_i)$ given by proposition
1.9 for the sequence $(\varepsilon_i)$.  

Choose an even integer $m\in \IN$ with
$$m > 8C/\delta\ .\leqno(2.2)$$

From the theory of spreading models there 
exists $(z_i)^{2m}_{i=1}$, a finite subsequence of $(y_i)$, such
that setting $\lambda = \Vert \sum^{2m}_{i=1} z_i\Vert^{-1}$,
$$2 > \lambda \Vert \sum_{i\in F} z_i\Vert > 1/3\ .\leqno(2.3)$$
whenever $F \subseteq \{1, \cdots, 2m\}$ with $\vert F\vert \ge m$.  

Thus there exists $x = \sum^{2m}_{i=1} x_i\in 2C Ba S$ with $(x_i)$ a
block basis of $(e_i)$ and 

\noindent $\Vert Tx_i - \lambda z_i\Vert <\varepsilon_i$
for $i \le 2m$.  For $i\le 2m$ choose $\omega_i\in S$ with $T\omega_i =
\lambda z_i -Tx_i$ and $\Vert \omega_i\Vert \le C \varepsilon_i$.  Hence
$T(x_i +\omega_i) = \lambda z_i$.

Since $\Vert T(\sum^{2m}_1 (x_i +\omega_i))\Vert > 1/3$, and
$$\Vert \sum^{2m}_1 (x_i +\omega_i)\Vert \le \Vert \sum^{2m}_1 x_i\Vert +
\sum^{2m}_1 \Vert \omega_i\Vert < 2C +\sum^\infty_1 \varepsilon_i C <3C\ ,$$
by lemma 2.2 we have $\Vert \sum^{2m}_1(x_i +\omega_i)\Vert_0 >\delta$.
Since $\Vert \sum^{2m}_1\omega_i\Vert_0 \le
\Vert \sum^{2m}_1\omega_i\Vert <\delta/2$ by (2.1) there exists
$i_1 \le 2m$ with $\Vert x_{i_1}\Vert_0 > \delta/2$.

Now 
$$\Vert T(\sum^{2m}_{\scriptstyle i=1\atop\scriptstyle  i\ne i_1}  
(x_i +\omega_i))\Vert =
\Vert \sum^{2m}_{\scriptstyle i=1\atop \scriptstyle i\ne i_1} 
\lambda z_i\Vert > {1\over 3}$$
and so we may repeat the argument above finding $i_2 \ne i_1$ with $\Vert
x_{i_2}\Vert_0 > \delta/2$.  In fact by (2.3) we can repeat this $m$-times
obtaining distinct integers $(i_k)^m_{k=1}  
\subseteq \{1,2,\cdots, 2m\}$ with
$\Vert x_{i_k}\Vert_0 >\delta/2$ for $k\le m$.  But then $2C\ge \Vert x\Vert
= \Vert \sum^{2m}_{i=1} x_i\Vert \ge \Vert 
\sum^m_{k=1} x_{i_k}\Vert \ge \sum^m_{k={m\over2}+1}
\Vert x_{i_k}\Vert_0 \ge
\delta m/4$ which contradicts (2.2).\qed
\bigskip

\noindent {\bf \S 3. Open Problems}
\medskip

Our work suggests a number of problems, of which we list a few.  For a more
extensive list of related problems and an overview of the current state of
infinite dimensional Banach space theory, see \rf{R}.

\noindent{\bf Problem 1}.  Let $X$ be a Banach space having property (WU) 
which does not contain $\ell_1$ and let $Y$ be a quotient of $X$.  Does
$Y$ have property (WU)?

In light of Theorem A it is worth noting that $C(\omega^\omega)$ has property
(WU) \rf{MR} but does not embed into any space having a shrinking unconditional
f.d.d.  In fact $C(\omega^\omega)$ is not even a subspace of a quotient of
such a space.  Indeed $C(\omega^\omega)$ fails property (U) (see {\it e.g.}
\rf{HOR}) while any quotient of a space with a shrinking unconditional f.d.d.
will have property (U).  In fact if $X$ has property (U) and does not 
contain $\ell_1$,
then any quotient of $X$ will have property (U) \rf{R}.  The next problem is
due to H.\ Rosenthal.
\medskip
\noindent{\bf Problem 2}.  Let $X$ have a shrinking unconditional f.d.d.
and let $Y$ be a quotient of $X$.  Does $Y$ embed into a Banach space
having a shrinking unconditional f.d.d.?

We say that a Banach space $Y$ has uniform-(WU) if there exists $K <\infty$
such that every normalized weakly null sequence in $Y$ has a $K$-unconditional
subsequence.  Our proof of Theorem A showed that the quotient space $Y$ has 
uniform-(WU).
\medskip
\noindent{\bf Problem 3}.  If $Y$ has property (WU) does $Y$ have uniform-(WU)?

Theorem B solved a special case of the following well known problem.
\medskip
\noindent{\bf Problem 4}.  Let $Y$ be a quotient of $C(\omega^\omega)$ (or
more generally $C(K)$ where $K$ is a compact countable metric space).  Is
$Y$ $c_0$-saturated?

Regarding this problem, T.\ Schlumprecht \rf{Sc} has observed that if $Y$ is 
a quotient of $C(\omega^\omega)$, then the closed linear span of any
normalized weakly null sequence in $Y$ which has $\ell_1$ as a spreading
model must contain $c_0$.

It is not true that the quotient of a $c_0$-saturated space must also be
$c_0$-saturated.  The separable Orlicz function space $H_M(0,1)$, with $M(x) =
(e^{x^4}-1)/(e-1)$, considered in \rf{CKT} is $c_0$-saturated and yet has 
$\ell_2$ as a quotient.  We wish to thank S.\ Montgomery-Smith for bringing
this fact to our attention.  However this space does not have an unconditional
basis and so we ask
\medskip
\noindent{\bf Problem 5.}  Let $X$ be a $c_0$-saturated space with an
unconditional basis and let $Y$ be a quotient of $X$.  Is $Y$ $c_0$-saturated?

A more restricted and perhaps more accessible question
is the following
\medskip
\noindent{\bf Problem 6}.  Let $Y$ be a quotient of $S_n$, the $n^{th}$-Schreier
space, where $n \ge 2$.  Is $Y$ $c_0$-saturated?
Does $Y$ have property-$S(n)$?

$S_n$ is defined as follows.  Let $\Vert x\Vert_1$ be the Schreier norm.
If $(S_n, \Vert \cdot \Vert_n)$ has been defined, set for $x\in c_{00}$,
the finitely supported real sequences,
$$\Vert x\Vert_{n+1} = \max \{\sum^p_{k=1} \Vert E_k x\Vert_n: p \le E_1
< E_1 < \cdots < E_p\}\ .$$
(Here $p \le E_1$ means $p\le \min E_1$ and $E_1 < E_2$ means $\max E_1 <
\min E_2$.  Also $Ex(i) = x(i)$ if $i\in E$ and 0 otherwise.)  $S_{n+1}$
is the completion of $(c_{00}, \Vert \cdot \Vert_{n+1})$.  
The unit vector basis $(e_n)$ is a
1-unconditional shrinking basis for every $S_n$ and $S_n$ embeds into
$C(\omega^{\omega^n})$.

Property-$S(n)$ is defined as follows.  {\it $n$-averages} of a sequence
$(y_m)$ are defined inductively:  an $n+1$-average of $(y_m)$ is a
1-average of a block basis of normalized $n$-averages.  $Y$ has 
{\it property-$S(n)$} if every normalized weakly null basic sequence in $Y$
admits a block basis of $n$-averages equivalent to the unit vector basis
of $c_0$.  $S_n$ has property-$S(n)$.
\bigskip

\baselineskip=12pt\frenchspacing
\centerline{\bf References}
\medskip
\item{[BL]}  B. Beauzamy and J.-T. Laprest\'e, 
{\it Mod\`eles \'etal\'es
des espaces de Banach}, Travaux en Cours, Hermann, Paris (1984).

\item{[CKT]}  P.G. Casazza, N.J. Kalton and L. Tzafriri, 
{\it Decompositions of Banach lattices into direct sums}, 
Trans. Amer. Math. Soc. {\bf 304} (1987), 771--800.

\item{[HOR]}  R. Haydon, E. Odell and H. Rosenthal, {\it On certain classes
of Baire-1 functions with applications to Banach space theory}, preprint.

\item{[J]}  W.B. Johnson, {\it On quotients of $L_p$ which are quotients
of $\ell_p$}, Compositio Math. {\bf 34} (1977), 69--89.

\item{[Ja]}  R.C. James, {\it Uniformly non-square Banach spaces}, Ann. of
Math. {\bf 80} (1964), 542--550.

\item{[JZ1]}  W.B. Johnson and M. Zippin, {\it On subspaces of quotients
of $(\sum Gn)_{\ell p}$ and $(\sum Gn)_{c_0}$}, Israel J. Math. {\bf 13} 
(1972), 311--316.

\item{[JZ2]}  W.B. Johnson and M. Zippin, {\it Subspaces and quotients
of $(\sum G_n)_{\ell p}$ and $(\sum G_n)_{c_0}$}, Israel J. Math. 
{\bf 17} (1974), 50--55.

\item{[LT1]}  J. Lindenstrauss and L. Tzafriri, 
``Classical Banach Spaces I,'' Springer-Verlag, Berlin, 1977.

\item{[LT2]}  J. Lindenstrauss and L. Tzafriri,
``Classical Banach Spaces II,'' Springer-Verlag, Berlin, 1979.

\item{[MR]}  B. Maurey and H. Rosenthal, {\it Normalized weakly null
sequences with no unconditional subsequence}, Studia Math. 
{\bf 61} (1977), 77--98.

\item{[O]}  E. Odell, {\it A normalized weakly null sequence with no
shrinking subsequence in a Banach space not containing $\ell_1$}, 
Composite Math. {\bf 41} (1980), 287--295.

\item{[R]}  H. Rosenthal, {\it Some aspects of the subspace structure of
infinite dimensional Banach spaces}, in ``Approximation Theory and
Functional Analysis'' (ed. C.~Chuy), Academic Press, 1991, 151--176.

\item{[S]}  J. Schreier, {\it Ein Gegenbespiel zur Theorie der schwachen
Konvergenz}, Studia Math. {\bf 2} (1930), 58--62.

\item{[Sc]}  T. Schlumprecht, private communication.

\item{[Z]}  M. Zippin, {\it Banach spaces with separable duals},
Trans. Amer. Math. Soc. {\bf 310} (1988), 371--379.

\end